\documentclass[12pt]{article}

\usepackage{graphicx}

\def\ra{\rightarrow}

\def\ss{\subseteq}

\def\d{\delta}

\def\O{\Omega}
\def\sss{\subset \! \subset}

\def\dbar{\overline{\partial}}

 \def\HollowBox #1#2{{\dimen0=#1 \advance\dimen0 by -#2       
       \dimen1=#1 \advance\dimen1 by #2                       
        \vrule height #1 depth #2 width #2                    
        \vrule height 0pt depth #2 width #1                   
        \llap{\vrule height #1 depth -\dimen0 width \dimen1}%
       \hskip -#2                                             
       \vrule height #1 depth #2 width #2}}                   
 \def\BoxOpTwo{\mathord{\HollowBox{6pt}{.4pt}}\;}             

\def\endpf{\hfill $\BoxOpTwo$}

\font\teneufm=eufm10
\font\seveneufm=eufm7
\font\fiveeufm=eufm5
\newfam\eufmfam
\textfont\eufmfam=\teneufm
\scriptfont\eufmfam=\seveneufm
\scriptscriptfont\eufmfam=\fiveeufm

\newfam\msbfam
\font\tenmsb=msbm10  scaled \magstep1 \textfont\msbfam=\tenmsb
\font\sevenmsb=msbm7 scaled \magstep1 \scriptfont\msbfam=\sevenmsb
\font\fivemsb=msbm5  scaled \magstep1 \scriptscriptfont\msbfam=\fivemsb
\def\Bbb{\fam\msbfam \tenmsb}

\def\RR{{\Bbb R}}
\def\CC{{\Bbb C}}

\def\NN{{\Bbb N}}

\newtheorem{theorem}{Theorem}
\newtheorem{corollary}[theorem]{Corollary}
\newtheorem{proposition}[theorem]{Proposition}
\newtheorem{lemma}[theorem]{Lemma}
\newtheorem{remark}[theorem]{Remark}

\newtheorem{definition}{Definition}
\newtheorem{example}[definition]{EXAMPLE}

\makeindex

\begin{document}

\begin{center}
\huge \bf Complex Analysis as Catalyst
\end{center}

\begin{center}
\large Steven G. Krantz\footnote{We are happy to thank
the American Institute of Mathematics for its hospitality
and support during this work.}
\end{center}
\vspace*{.15in}

\begin{quote}
{\bf Abstract:}  We see the subject of complex
analysis, in both one and several variables, as
an excuse to study other areas:  real variables, topology,
group theory, abstract algebra, partial differential equations, or geometry.
The purpose of this paper is to develop that theme, illustrated
by several examples.
\end{quote}

\setcounter{section}{-1}

\section{Prefatory Thoughts}

I am known as a complex analyst, yet I have
never had a course in complex variable theory.  This lacuna
in my education has perhaps had an effect on the way that
I view the subject.

I see the elegant subject of complex function theory as
a stage for the setting of beautiful problems.
Once these problems are apprehended or understood, they
are often best solved by stripping away the complex variables
and bringing in tools from other parts of mathematics. 
That is the point of view that we promulgate in the present 
article.  We illustrate the point with examples from
real analysis, group theory, abstract algebra, partial differential equations,
geometry, and other subjects as well.

\section{An Example from Topology}

Topology has become part of the bedrock of modern mathematics.
Originally studied as {\it analysis situs}, it now permeates
all branches of our discipline.  We understand that the
entire texture of a subject may be influenced by the particular
topology that is chosen.  The examples in this section illustrate
the point.

In this paper a {\it domain} will be a connected, open set in
$\CC$ or $\CC^n$. Let $\Omega$ be a domain. It is a
commonplace in complex variable theory to introduce the notion
of {\it uniform convergence on compact subsets} of $\Omega$.
Though not often noted, this is in fact equivalent to the
compact-open topology on the space ${\cal O}(\Omega)$ of
holomorphic functions on $\Omega$.  In practice we say
that $f_j \in {\cal O}(\Omega)$ converge uniformly on compact
subsets of $\Omega$ to a function $f \in {\cal O}(\Omega)$ if,
for any compact $K \subseteq \Omega$ and any $\epsilon > 0$, there
is a $J >0$ such that if $j > J$ then $|f_j(z) - f(z)| < \epsilon$
for all $z \in K$.  The basic result, which every student learns,
is that the limit function $f$ is holomorphic on $\Omega$ provided
only that the elements $f_j$ are holomorphic on $\Omega$.

An interesting and natural question to ask at this point is
\begin{quote}
What happens if the $f_j \in {\cal O}(\Omega)$ only
converge pointwise?
\end{quote}
A naive answer might be ``Nothing much,'' but that would
be hasty and also incorrect.  Certainly the limit function
$f$ is well defined and we may ask what properties it has.
Is it continuous?  Probably not.  It turns out that
the correct statement is this:

\begin{theorem} \sl
Let $f_j \in {\cal O}(\Omega)$ converge pointwise to a
limit function $f$ on $\Omega$.  Then there is a dense,
open subset $V \subseteq \Omega$ such that $f$ is holomorphic
on $V$.
\end{theorem}

As one might suppose from the statement, the proof will involve
the Baire category theorem.  We see then, in our first example,
that even though the theorem is about function theory
the proof will consist of real analysis and topology.
\medskip \\

\noindent {\bf Proof of Theorem 1:}  Let $U$ be an open subset
of $\Omega$ with compact closure in $\Omega$.  Define, for 
$k = 1, 2, \dots$,
$$
S_k = \{z \in \overline{U}: |f_j(z)| \leq k \ \hbox{for all} \ j \in \NN\} \, .
$$
Since the $f_j$ converge at each $z \in \overline{U}$, certainly the set
$\{f_j(z)\}$ is bounded for each fixed $z$.  
So each $z \in \overline{U}$ lies in some $S_k$.
In other words,
$$
\overline{U} = \bigcup_k S_k \, .
$$

Now of course $\overline{U}$ is a complete metric space (in the
ordinary Euclidean metric), so the Baire category theorem tells
us that some $S_k$ must be ``somewhere dense'' in $\overline{U}$.
This means that $\overline{S_k}$ will contain a nontrivial Euclidean metric
ball in $\overline{U}$.  Call the ball $B$.  Now it is a simple
matter to apply Montel's theorem on $B$ to find a subsequence
$f_{j_k}$ that converges uniformly on compact sets to a limit
function $g$.  But of course $g$ must coincide with $f$, and
$g$ (hence $f$) must be holomorphic on $B$.  

Since the choice of $U$ in the above arguments was arbitrary, the conclusion
of the theorem follows.
\endpf
\medskip \\

One might ask what more can be said about the open, dense set on
which the limit function $f$ is holomorphic.  Put in other words,
what can one say about $\Omega \setminus V$?  In fact Siciak [SIC]
has given a rather thorough answer to this question.  The statement
and proofs involve analytic capacity, and we cannot treat the details here.

\section{Hard Analysis}

The central feature of the basic theory of complex analysis is
the Cauchy integral formula.  A basic form of this result is
as follows:

\begin{proposition} \sl
Let $f$ be holomorphic in a neighborhood of
the closed disc $\overline{D}(P,r)$ in the complex
plane.  If $z \in D(P,r)$ then
$$
f(z) = \frac{1}{2\pi i} \oint_{\partial D(P,r)} \frac{f(\zeta)}{\zeta - z} \, d\zeta \, .  \eqno (*)
$$
\end{proposition}

A reasonable---though not often formulated---question to ask is,
``What is so special about holomorphic functions?  Is the formula
$(*)$ true only for holomorphic functions?  Is something analogous
true for a more general class of functions?  Is there a universal
formula of which $(*)$ is just a special case?''

The somewhat surprising answer to these questions is ``yes''.  We shall,
in the present section, prove that result.  Later sections will give
applications of the more general integral formula.

The experienced reader will know that the Poisson integral formula is
an identity that is closely related to $(*)$ (see [GRK], [KRA4], [KRA5]), and which is valid for
harmonic functions.  So it can properly be considered a generalization
of the Cauchy formula.  But that is not what we seek at this time.  First
of all, it is too well known hence not suitable grist for our mill.  Secondly,
we seek a formula that is valid for all smooth functions.  In fact our
result is this:

\begin{theorem} \sl.
Let $\Omega \subseteq \CC$ be a bounded domain with $C^1$ boundary (i.e.,
the boundary consists of finitely many closed, $C^1$ curves).  Let
$f$ be a continuously differentiable function on $\overline{\Omega}$.  Then,
for any $z \in \Omega$,
$$
f(z) = \frac{1}{2\pi i} \oint_{\partial \Omega} \frac{f(\zeta)}{\zeta - z} \, d\zeta -
   \frac{1}{2\pi i} \mathop{\int \!\!\!\int}_\Omega \left ( 
             \frac{\partial f/\partial \overline{\zeta}}{\zeta - z} \right ) \, d\overline{\zeta} \wedge d\zeta \, .
$$
[In the proof we shall review the notation $\partial\partial \overline{\zeta}$.]
\end{theorem}

Our proof will use a form of Stokes's theorem that may be unfamiliar.  Recall that
the classical Stokes's theorem that can be found in any calculus book (see, for example,
[BLK]) says that, if $u_1$, $u_2$ are continuously differentiable functions on $\overline{\Omega}$,
then
$$
\int_{\partial \Omega} u_1 dx +  \int_{\partial \Omega} u_2 dy =
   \mathop{\int \!\!\! \int}_\Omega \frac{\partial u_2}{\partial x} - \frac{\partial u_1}{\partial y} \ dx \wedge dy \, . \eqno (\star)
$$
Note that, in complex notation, 
$$
\frac{\partial}{\partial z} = \frac{1}{2} \left ( \frac{\partial}{\partial x} - i \frac{\partial}{\partial y} \right ) \, ,
$$
$$				   
\frac{\partial}{\partial \overline{z}} = \frac{1}{2} \left ( \frac{\partial}{\partial x} + i \frac{\partial}{\partial y} \right ) \, ,
$$
$dz = dx + i dy$, and $d\overline{z} = dx - i dy$.  Thus a simple change of variable in the
formula $(\star)$ yields the ``complex'' version of Stokes's theorem
$$
\oint_{\partial \Omega} u dz = \mathop{\int \!\!\!\ \int}_\Omega \frac{\partial u}{\partial \overline{z}} \ d\overline{z} \wedge dz \, .
\eqno (\star\star)
$$
We shall make good use of this version of Stokes's result in the argument that follows.
\medskip \\

\noindent {\bf Proof of Theorem 3:}  We shall unashamedly use Stokes's theorem, and a mild dose of the theory of
differential forms.  This treatment will be self-contained, and even the nervous
reader should find it enlightening.

Fix a point $z \in \Omega$ and a positive number $\epsilon$ that is less than the Euclidean
distance of $z$ to $\partial \Omega$.  Define
$$
\Omega_\epsilon \equiv \Omega \setminus \overline{D}(z, \epsilon) \, .
$$
We apply the aforementioned version of Stokes's theorem to the function
$$
g(\zeta) \equiv \frac{f(\zeta)}{\zeta - z} 
$$
on the domain $\Omega_\epsilon$.

The result is
$$
\oint_{\partial \Omega_\epsilon} g(\zeta) \, d\zeta =
   \mathop{\int \!\!\!\int}_{\Omega_\epsilon} \frac{\partial}{\partial \overline{\zeta}} \left (
        g(\zeta ) \right ) \, d\overline{\zeta} \wedge d\zeta \, .
$$
Writing out the definition of $g$ gives
$$ 
\oint_{\partial \Omega_\epsilon} \frac{f(\zeta)}{\zeta - z} \, d\zeta =
   \mathop{\int \!\!\! \int}_{\Omega_\epsilon} \frac{\partial}{\partial \overline{\zeta}} \left (
      \frac{f(\zeta)}{\zeta - z} \right ) d\overline{\zeta} \wedge d\zeta
$$
or
$$
\oint_{\partial \Omega_\epsilon} \frac{f(\zeta)}{\zeta - z} \, d\zeta =
   \mathop{\int \!\!\! \int}_{\Omega_\epsilon} \left (
      \frac{\partial f/\partial \overline{\zeta}}{\zeta - z} \right ) d\overline{\zeta} \wedge d\zeta \, .  \eqno (\dagger)
$$

We may write out the lefthand side of $(\dagger)$ as
$$
\oint_{\partial \Omega} \frac{f(\zeta)}{\zeta - z} \, d\zeta - 
   \oint_{\partial D(z, \epsilon)} \frac{f(\zeta)}{\zeta - z} \, d\zeta \equiv I + II \, .
$$
Note the differing orientations of the two pieces of the boundary.
We examine $II$:
\begin{eqnarray*}
II & = & \int_0^{2\pi} \frac{f(z + \epsilon e^{it})}{\epsilon e^{it}} \, \cdot \, i \epsilon e^{it} \, dt \\
   & = & i \int_0^{2\pi} f(z + \epsilon e^{it})  \, dt \, .
\end{eqnarray*}
A simple limiting argument, using the continuity of $f$, shows that this last line tends as $\epsilon \ra 0^+$
to $2\pi i f(z)$.  

Thus we have from $(\dagger)$ (letting $\epsilon \ra 0^+$ on the right as well) that
$$
\oint_{\partial \Omega} \frac{f(\zeta)}{\zeta - z} \, d\zeta - 2\pi i f(z) =
    \mathop{\int\!\!\!\int}_\Omega \left ( \frac{\partial f/\partial \overline{\zeta}}{\zeta - z} \right ) 
    \, d\overline{\zeta} \wedge d\zeta \, .
$$
In conclusion,
$$
f(z) = \frac{1}{2\pi i} \oint_{\partial \Omega} \frac{f(\zeta)}{\zeta - z} \, d\zeta -
   \frac{1}{2\pi i} \mathop{\int \!\!\! \int}_\Omega \left ( \frac{\partial f/\partial \overline{\zeta}}{\zeta - z} 
   \right )  \,  d\overline{\zeta} \wedge d\zeta \, .
$$
That is the desired conclusion.
\endpf 
\medskip \\

Of course it is incumbent on us to point out that, in the case that $f$ is holomorphic, it follows
immediately (from the Cauchy-Riemann equations) that $\partial f/\partial \overline{\zeta} \equiv 0$.
Hence the formula in the statement of the theorem becomes
$$
f(z) = \frac{1}{2\pi i} \oint_{\partial \Omega} \frac{f(\zeta)}{\zeta - z} \, d\zeta  \, .
$$
This is of course the classical Cauchy integral formula.

\section{Partial Differential Equations}

One of the most profound facts in basic function theory is the 
Riemann mapping theorem:

\begin{theorem} \sl
Let $\Omega$ be any simply connected domain in $\CC$ that is not the
entire complex plane.  Then there is a one-to-one, onto conformal 
mapping of $\Omega$ to the unit disc $D$.
\end{theorem}

This result has exerted a profound influence on our subject ever since
its original formulation (and not-quite-right) proof.  It has been
generalized to Koebe's uniformization theorem, and has inspired the
definition of the Carath\'{e}odory and Kobayashi metrics (see, for instance,
[KRA1] and also our Section 7).  If $\varphi: \Omega \ra D$ is the Riemann mapping, then it
is natural to think that the function theory on $\Omega$ may be transferred
to function theory on $D$ by way of $\varphi$.  This process has many advantages,
for $D$ has great symmetry, it is possessed of a transitive automorphism group,
and its boundary carries the full arsenal of Fourier analysis.

But many of the most basic questions that one might wish to ask would necessitate
that $\varphi$ be well behaved up to the boundary of $\Omega$.  Thus one might
ask whether $\varphi$ and $\varphi^{-1}$ extend continuously to $\partial \Omega$
and $\partial D$ respectively.  Carath\'{e}odory has shown (see the details in [GRK])
that, in case $\partial \Omega$ is a Jordan curve, the answer is ``yes''.  More generally,
one might ask whether, if $\partial \Omega$ is $C^k$ (i.e., consists of the union
of finitely many $C^k$ curves), then do $\varphi$ and $\varphi^{-1}$ extend
$C^k$ to the respective boundaries?   Phrased slightly differently, we may ask whether
$\varphi$ and $\varphi^{-1}$ extend to $C^k$ diffeomorphisms of the closures of
$\Omega$ and $D$.

These last are profound questions, first studied in the nineteenth century 
by Painlev\'{e} in his Paris
thesis.  Painlev\'{e}'s original techniques were purely function-theoretic, but
the modern approach to the matter is by way of partial differential equations.
A full treatment of the matter requires extensive exposition (see the
monograph [KRA2]), but we may indicate the main points here.

In fact the approach that we wish to present requires a detour into the study
of the celebrated Dirichlet problem.  Now let $\Omega \subseteq \CC$ be a bounded
domain with $C^1$ boundary.  Let $f$ be a given continuous function on $\partial \Omega$.
The {\it Dirichlet problem} is the elliptic boundary value problem given by
$$
\left \{ \begin{array}{lcl}
    \bigtriangleup u & = & \quad 0 \qquad \qquad \hbox{on} \quad \Omega \\
                   u & = & \quad f \qquad \qquad \hbox{on} \quad \partial \Omega \, .
	 \end{array}
\right.
$$
As a classical instance, if $\Omega$ is a metal sheet of heat-conducting material,
and if $f$ represents an initial heat distribution on the boundary of $\Omega$,
then the solution $u$ of the Dirichlet problem is the steady state heat distribution
over $\Omega$.  

There is an extensive theory---due to O. Perr\`{o}n (see [GRK])---for constructing the solution
of the Dirichlet problem.  In the case that $\Omega$ is the unit disc, the Poisson integral
formula gives an explicit solution to the Dirichlet problem:
$$
u(r e^{i\theta}) = \frac{1}{2\pi} \int_0^{2\pi} f(e^{i\psi}) \frac{1 - r^2}{1 - 2r \cos(\theta - \psi) + r^2} \, d\psi \, .
$$
One notes immediately in this formula that the kernel, that is to say the expression
$$
P(r, \theta- \psi) \equiv \frac{1 - r^2}{1 - 2r \cos(\theta - \psi) + r^2} \, ,
$$
is real analytic in its arguments.  It follows then that any harmonic function is real
analytic, so certainly smooth.  The key question is
\begin{quote}
What is the behavior of the solution $u$ of the Dirichlet problem
{\it up to the boundary} of $\Omega$?
\end{quote}

Perr\`{o}n's original theory shows that, under mild regularity conditions on $\partial \Omega$,
the solution to the Dirichlet problem is continuous on the closure.  That is to say,
the function
$$
U(z) = \left \{ \begin{array}{lcl}
               u(z) & \hbox{if} & \quad z \in \Omega \\
	       f(z) & \hbox{if} & \quad z \in \partial \Omega
		\end{array}
       \right.
$$
is continuous on $\overline{\Omega}$. This assertion may be proved by
direct estimation from the Poisson integral formula in the case that
$\Omega$ is the disc (see [KRA3, Chapter 8]). The higher (or $C^k$)
regularity theory---commonly known as the {\it Schauder estimates}---is
deep and difficult. It is challenging even to find a source that exposits
the matter completely and clearly.\footnote{The book [KRA2] is one such reference;
see also the classical references [HOR] and [GIT].}  We can only content
ourselves here with a brief statement (no proof) of the regularity theory
for the Laplacian.

First some terminology.  Recall that if $g$ is a function on a domain
$E \subseteq \RR^N$ then we say that $g$ satisfies a {\it classical
Lipschitz condition} if there is a constant $C > 0$ such that
$$
|g(x) - g(y)| \leq C \cdot |x - y|
$$
for all $x, y \in E$.  We sometimes call this a Lipschitz condition of {\it order 1}.
More generally, if $0 < \alpha < 1$, then 
we say that $g$ satisfies a {\it Lipschitz condition
of order $\alpha$} provided that
$$
|g(x) - g(y)| \leq C \cdot |x - y|^\alpha
$$
for all $x, y \in E$.\footnote{It seems to be the
case that harmonic analysts like to say ``Lipschitz space'' while partial differential
equations theorists like to say ``H\"{o}lder space''.  We adhere to the former
paradigm.}  Clearly a Lipschitz condition of order 1 is stronger than
a Lipschitz condition of order $\alpha$, $0 < \alpha < 1$.  

Now let $g$ be a function on an open domain $U \subseteq \RR^N$. Let $k$ be
a nonnegative integer and $0 < \alpha < 1$. We say that $g \in
C^{k,\alpha}$ provided that $g$ is $k$-times continuously differentiable
on $U$ and further that any $k^{\rm th}$ derivative $D^\beta g$ (with
$\beta$ a multi-index of order $k$) is Lipschitz\footnote{This definition is easily
extended to the situation where the domain of the function is the closure of
an open domain.}  of order $\alpha$.  The $C^{k,\alpha}$
form a scale of spaces that stratify $C^\infty$.  Clearly, for example, for fixed $\alpha$,
$$
\bigcap_k C^{k, \alpha} = C^\infty \, .
$$
Now we may state our main regularity result for the Dirichlet problem for the Laplacian.
Again, references are [KRA2], [GIT], and [HOR].

\begin{proposition} \sl
Let $k$ be a nonnegative integer.  Let $\Omega \ss \CC$ be a bounded domain with $C^k$ boundary.
Let $f$ be a $C^k$ function on $\partial \Omega$.  Then the solution $u$ on $\Omega$ of
the Dirichlet problem with boundary data $f$ satisfies $u \in C^{(k-1),\alpha}(\overline{\Omega})$
for any $0 < \alpha < 1$.
\end{proposition}

The fact that we may not assert in this last proposition that $u \in C^k(\overline{\Omega})$ is
an unfortunate artifact of the theory of singular integrals (again see [KRA2] and also [STE]).
A sharper result may be formulated if we use Sobolev spaces or Besov spaces or Triebel-Lizorkin
spaces.  One may even formulate a more refined definition of Lipschitz spaces that results
in a sharper statement.  But we cannot indulge in these technical niceties here.

We next need the following sharp form of the classical 
lemma of Hopf (which was originally formulated to
study the maximum principle for solutions of elliptic partial
differential equations):

\begin{lemma}[Hopf]  \sl
Let $\Omega \sss \RR^N$ have $C^2$ boundary.  Let $u \in C(\bar \Omega)$
with $u$ harmonic and non-constant on $\Omega.$  Let $P \in \bar \Omega$ and assume that
$u$ takes a local minimum at $P.$  Then
$$
\frac{\partial u}{\partial \nu} (P) < 0 . 
$$
Here $\partial/\partial \nu$ is the unit outward normal derivative.
\end{lemma}

\noindent {\bf Remark:}  That the indicated normal derivative is nonpositive
is obvious just from examining the Newton quotients.  The interesting
thing is that the derivative must be strictly negative.
\endpf
\medskip \\

\noindent {\bf Proof:}  Suppose without loss of generality that
$u > 0$ on $\Omega$ near $P$ and that 
$u(P) = 0.$  Let $B_R$ be a ball that is internally
tangent to $\bar \Omega$ at $P.$  We may assume that
the center of this ball is at the origin and that
$P$ has coordinates $(R,0,\dots,0).$
Then, by Harnack's inequality
(see [GRK]), we have for $0 < r < R$ that
$$
u(r,0,\dots,0) \geq c \cdot \frac{R^2 - r^2}{R^2 + r^2}
$$
hence
$$
\frac{u(r,0,\dots,0) - u(R,0,\dots,0)}{r - R} \leq - c' < 0 .
$$
Therefore
$$
\frac{\partial u}{\partial \nu}(P) \leq - c' < 0 . 
$$
This is the desired result.
\endpf 
\smallskip \\

\begin{proposition} \sl
Let $\Omega \subseteq \CC$ be a bounded, simply connected domain with $C^k$ boundary.
Let $\varphi: \Omega \ra D$ be the conformal mapping provided by the
Riemann mapping theorem.  Then $\varphi$ extends to a $C^{k-2}$, univalent,
mapping of $\overline{\Omega}$ to $\overline{D}$ for any $0 < \alpha < 1$.
\end{proposition}
{\bf Proof:}  Let $W$ be a collared neighborhood of $\partial \Omega$ with
smooth boundary.  Set $\partial \Omega' = \partial W \cap \Omega$ and let $\partial D'= \varphi(\partial \Omega').$
Define $B$ to be the region bounded by $\partial D$ and $\partial D'.$
Let $U$ be the region bounded by $\partial \Omega$ and $\partial \Omega'$.
We solve the Dirichlet problem on $B$ with boundary data
$$
f(\zeta) = \left \{ \begin{array}{lcl}
                       1 & \mbox{if} & \zeta \in \partial D \\
                       0 & \mbox{if} & \zeta \in \partial D'
                    \end{array}
           \right. 
$$
Call the solution $u.$  

Consider $v \equiv u \circ \varphi: \Omega \ra \RR.$  Then of course
$v$ is still harmonic.  By Carath\'{e}odory's theorem,
$v$ extends to $\partial \Omega, \partial \Omega',$ and
$$
v = \left \{ \begin{array}{lcl}
                       1 & \mbox{if} & \zeta \in \partial \Omega \\
                       0 & \mbox{if} & \zeta \in \partial \Omega' \\
                    \end{array}
           \right. 
$$
The function $v$ will be $C^{(k-1), \alpha}$, $0 < \alpha < 1$, on $\overline{U}$ by Proposition 5.
If we
consider a first order derivative ${\cal D}$ of $v$ we obtain
$$
|{\cal D} v | = |{\cal D}(u \circ \varphi)| = |\nabla u| \, |\nabla \varphi| \leq C . 
$$
It follows that
$$
|\nabla \varphi | \leq \frac{C}{|\nabla u|} . \eqno (*) 
$$

Now let us return to the $u$ from the Dirichlet problem that we considered
prior to line $(*)$.  Hopf's lemma tells us that $|\nabla u | \geq c' > 0$
near $\partial D.$  Thus, from $(*)$, we conclude that
$$
|\nabla \varphi| \leq C . \eqno (**) 
$$
Thus we have bounds on the first derivatives of $\varphi.$  

To control the second derivatives, we calculate that
\begin{eqnarray*}
C & \geq & |\nabla^2 v| = |\nabla(\nabla v)| = |\nabla(\nabla(u \circ \varphi))| \\
  & = &    |\nabla(\nabla u(\varphi) \cdot \nabla \varphi) | 
     = |\bigl (\nabla^2 u \cdot [\nabla \varphi]^2 \bigr ) 
            + \bigl (\nabla u \cdot \nabla^2 \varphi \bigr ) | . 
\end{eqnarray*}
Here the reader should think of $\nabla$ as representing a generic first
derivative and $\nabla^2$ a generic second derivative.
We conclude that
$$
|\nabla u| \, |\nabla^2 \varphi| \leq C + | \nabla^2 u | \, |(\nabla \varphi)^2 | \leq C' . $$
Hence (again using Hopf's lemma)
$$
|\nabla^2 \varphi| \leq \frac{C}{|\nabla u|} \leq C''.
$$ 

\noindent In the same fashion, we may prove that $|\nabla^j \varphi| \leq C_j,$ any
$0 \leq j \leq k - 1$  This means (use the fundamental theorem of calculus)
that $\varphi \in C^{k-2}(\bar \Omega).$
\endpf
\medskip \\

We note that several small technical emendations in the proof presented would
suffice to show that $\varphi \in C^{(k-1), \alpha}(\overline{\Omega})$ for
any $0 < \alpha < 1$.  We omit those details.

\section{Algebra}  

Now we turn to a result of Lipman Bers that has a different
flavor.  It characterizes domains in terms of a certain
algebraic invariant.

Let $\O \ss \CC$ be a domain. Let ${\cal O}(\O)$ denote the algebra of
holomorphic functions from $\O$ to $\CC$. Bers's theorem says, in effect,
that
the algebraic structure of ${\cal O}(\O)$ characterizes $\O$. We
begin our study by introducing a little terminology.

\begin{definition} \rm
Let $\O \ss \CC$ be a domain.  A $\CC$-algebra homomorphism
$\varphi: {\cal O}(\O) \ra \CC$ is called a {\it character}
of ${\cal O}(\O)$.  If $c \in \CC$, then the mapping
\begin{eqnarray*}
e_c: {\cal O}(\O) & \ra & \CC \, , \\
	       f  & \mapsto & f(c) \, , \\
\end{eqnarray*}
is
called a {\it point evaluation}.  Every point evaluation
is a character.
\end{definition}

It should be noted that if $\Omega$, $\widehat{\Omega}$ are domains and $\varphi: {\cal O}(\O) \ra {\cal O}(\widehat{\O})$ 
is not the trivial zero homomorphism, then $\varphi(1) = 1$.
This follows because $\varphi(1) = \varphi(1 \cdot 1) = \varphi(1) \cdot \varphi(1)$.
On any open set where the holomorphic function $\varphi(1)$ does not
vanish, we find that $\varphi(1) \equiv 1$.  The result follows
by analytic continuation.

It turns out that every character of ${\cal O}(\O)$ is a point
evaluation.  That is the content of the next lemma.

\begin{lemma} \sl
Let $\varphi$ be a character on ${\cal O}(\O)$.  Then 
$\varphi = e_c$ for some $c \in \O$.  Indeed, $c = \varphi(\hbox{id}) \in \O$.
Here $\hbox{id}$ is defined by $\hbox{id}(z) = z$.
\end{lemma}
{\bf Proof:} \  \ \ Let $c$ be defined as in the statement
of the lemma.  Let $f(z) = z - c$.  Then
$$
\varphi(f)  =  \varphi(\hbox{id}) - \varphi(c) 
	    =  c - c 
	    =  0 \, .
$$

If it were not the case that $c \in \O$ then the function
$f$ would be a unit in ${\cal O}(\O)$.  But then
$$
1 = \varphi(f \cdot f^{-1}) = \varphi(f) \cdot \varphi(f^{-1})
     = 0 \, .
$$
That is a contradiction.  So $c \in \O$.

Now let $g \in {\cal O}(\O)$ be arbitrary.  Then we may write
$$
g(z) = g(c) + f(z) \cdot \widetilde{g}(z) \, ,
$$
where $\widetilde{g} \in {\cal O}(\O)$.  Thus
$$
\varphi(g) = \varphi(g(c)) + \varphi(f) \cdot \varphi(\widetilde{g}) 
	   = g(c) + 0 = g(c) = e_c(g) \, .
$$
We conclude that $\varphi = e_c$, as was claimed.
\endpf
\medskip \\

Now we may state and prove Bers's theorem.

\begin{theorem} \sl
Let $\O$, $\widetilde{\O}$ be domains.  Suppose that
$$
\varphi: {\cal O}(\O) \ra {\cal O}(\widetilde{\O})
$$
is a $\CC$-algebra homomorphism.  Then there exists
one and only one holomorphic mapping $h: \widetilde{\O} \ra \O$
such that 
$$
\varphi(f) = f \circ h \quad \hbox{for all} \ f \in {\cal O}(\O) \, .
$$
In fact, the mapping $h$ is given by $h = \varphi(\hbox{id})$.

The homomorphism $\varphi$ is bijective if and only if $h$
is conformal, that is, a one-to-one and onto holomorphic mapping
from $\widetilde{\O}$ to $\O$.
\end{theorem}
{\bf Proof:}  Since we want the mapping $h$ to satisfy
$\varphi(f) = f \circ h$ for all $f \in {\cal O}(\O)$, it must
in particular satisfy $\varphi(\hbox{id}_\O) = \hbox{id}_\O \circ h = h$.
We take this as our definition of the mapping $h$.

If $a \in \widetilde{\O}$, then $e_a \circ \varphi$ is a character
of ${\cal O}(\O)$.  Thus our lemma tells us that
$e_a \circ \varphi$ must in fact be a point evaluation on $\O$.
As a result,
$$
e_a \circ \varphi = e_c \, , \quad \hbox{with} \ 
    c = (e_a \circ \varphi)(\hbox{id}_\Omega) = e_a(h) = h(a) \, .
$$
Thus, if $f \in {\cal O}(\O)$, then
$$
\varphi(f)(a) = e_a(\varphi \circ f) = (e_a \circ \varphi)(f)
     = e_{h(a)}(f) = f(h(a)) = (f \circ h)(a) 
$$
for all $a \in \widetilde{\O}$.
We conclude that $\varphi(f) = f \circ h$ for all $f \in {\cal O}(\O)$.  

For the last statement of the theorem, suppose that $h$ is a 
one-to-one, onto conformal mapping of $\widetilde{\O}$ to $\O$.
If $g \in {\cal O}(\O)$, then set $f = g \circ h^{-1}$.  It follows
that $\varphi(f) = f \circ h = g$.  Hence $\varphi$ is onto.
Likewise, if $\varphi(f_1) = \varphi(f_2)$, then $f_1 \circ h = f_2 \circ h$
hence, composing with $h^{-1}$, $f_1 \equiv f_2$.  So $\varphi$ is
one-to-one.  Conversely, suppose that $\varphi$ is an isomorphism.
Let $a \in \O$ be arbitrary.  Then $e_a$ is a character on ${\cal O}(\O)$;
hence $e_a \circ \varphi^{-1}$ is a character on ${\cal O}(\widetilde{\O})$.
By the lemma, there is a point $c \in \widetilde{\O}$ such
that $e_a \circ \varphi^{-1} = e_c$.  It follows that
$$
e_a = e_c \circ \varphi \, .
$$
Applying both sides of this last identity to $\hbox{id}_\O$ yields
$$
e_a(\hbox{id}_\O) = (e_c \circ \varphi)(\hbox{id}_\O) \, .
$$
Unraveling the definitions gives
$$
a \equiv e_c(\hbox{id}_\O \circ h) = h(c) \, .
$$
Thus $h(c) = a$ and $h$ is surjective.  The argument in fact
shows that the pre-image $c$ is uniquely determined.  So $h$
is also one-to-one.
\endpf 
\medskip \\

As an application\footnote{Of course this application is in a certain
sense trivial because the two domains have different topology.  All we are
doing here is illustrating the algebraic ideas with a different proof.}
of Bers's theorem, we can see immediately
that the disk $D = \{z \in \CC: |z| < 1\}$ and the annulus
${\cal A} = \{z \in \CC: 1/2 < |z| < 2\}$ are not conformally
equivalent.  For the algebra ${\cal O}(D)$ can be generated
by $1$ and $z$.  But the algebra ${\cal O}({\cal A})$ cannot
be generated by $1$ and just one other function (because
natural generators for ${\cal O}({\cal A})$ are $\{1, z, 1/z\}$ 
and it is impossible to come up with a shorter list).  We leave
the details of these assertions to the reader.

It must be noted, of course, that if $\Omega_1$ and $\Omega_2$ are
conformally equivalent via a mapping $\varphi: \Omega_1 \ra \Omega_2$,
then their respective algebras of holomorphic functions are trivially
isomorphic by way of the conjugation
$$
{\cal O}(\Omega_1) \ni \eta \longmapsto \varphi \circ \eta \circ \varphi^{-1} \in {\cal O}(\Omega_2) \, .
$$
The statement of the theorem is the converse of this trivial result, and
is definitely more interesting.  As an illustration, let $\Omega_1 = D$, the unit disc,
and $\Omega_2 = \CC$, the complex plane.  Of course these domains cannot be
conformally equivalent---by a simple application of Liouville's theorem.  We may conclude
from Theorem 9 that ${\cal O}(D)$ and ${\cal O}(\CC)$ are inequivalent as algebras.
This last statement is rather interesting, and is certainly nonobvious.  Both
algebras can be thought of as the set of all power series
$$
\sum_j a_j z^j \, ,
$$
and the only difference between the two is the radius of convergence.  It is not
so clear why the two algebras should be nonisomorphic.

\section{Group Theory}  

We now pass to the subject of the function theory of several complex variables.
One of the great classical results was motivated by the Riemann mapping
theorem from one complex variable.  That venerable result tells us that,
in a certain sense, the disc is the ``canonical domain'' in $\CC$.  Any
simply connected domain is conformally equivalent to the unit disc.
Thus one may ask what is the analogous result in several complex variables.
What is the canonical domain in $\CC^n$?

Two natural candidates for this canonical domain are the unit ball
$$
B = \{z = (z_1, z_2, \dots, z_n) \in \CC^n: \sum_j |z_j|^2 < 1\}
$$
and the unit polydisc
$$
D^n = \{z = (z_1, z_2, \dots, z_n) \in \CC^n: |z_j| < 1 \ \ \hbox{for} \ \ 
j = 1, 2, \dots, n\} \, .
$$
Poincar\'{e} caused quite a stir in 1906 when he proved that,
when $n > 1$, $B$ and
$D^n$ are {\it not} biholomorphically equivalent.  That is to say, there
does {\it not} exist a holomorphic mapping
$$
\Phi: B \ra D^n
$$
which is one-to-one and onto.  This result calls into question whether
there {\it could be} a canonical domain in $\CC^n$, $n > 1$.

We shall reproduce here Poincar\'{e}'s original proof. It is remarkable in
that it uses the automorphism group of a domain as a biholomorphic
invariant.

Let us recall now a few of the elementary ideas from the function theory of
several complex variables. Let $\Omega \ss \CC^n$ be a domain. A function
$f: \Omega \ra \CC$ is said to be {\it holomorphic} if it is holomorphic
in each variable separately. That is to say, if the values of $z_1, z_2,
\dots, z_{j-1}, z_{j+1}, \dots, z_n$ are frozen, then we mandate that
$\zeta \mapsto f(z_1, z_2, \dots, z_{j-1}, \zeta, z_{j+1}, \dots, z_n)$
be holomorphic as a function of one complex variable (whereever this
function makes sense). It can be shown, thanks to Hartogs (see [KRA3])
that such a function is $C^\infty$ as a function of several variables,
indeed it is real analytic as a multi-variable function. Certainly such a
function satisfies the Cauchy-Riemann equations.

Now if $\Omega$ is a domain then we let the {\it automorphism group} of $\Omega$,
denoted $\hbox{Aut}(\Omega)$, be the collection of holomorphic maps
$$
\Phi(z) = (\varphi_1(z), \dots, \varphi_n(z)) = (\varphi_1(z_1, \dots z_n), \dots,
\varphi_n(z_1, \dots, z_n))
$$
from $\Omega$ to itself which are both one-to-one and onto.  The collection
$\hbox{Aut}(\Omega)$ of such mappings forms a group when equipped with
the binary operation of composition of mappings.  We use the topology of
{\it uniform convergence on compact sets} (equivalently, the compact-open topology)
on the automorphism group.  Thus $\hbox{Aut}(\Omega)$ is a topological group, and
in fact (in case $\Omega$ is bounded), it can be shown to be a Lie group (see [KOB]).

It is important to note that the automorphism group is a biholomorphic invariant.  This
means the following.  Let 
$$
\Phi: \Omega_1 \ra \Omega_2
$$
be a biholomorphic mapping.  Then the induced mapping
$$
\hbox{Aut}(\Omega_1) \ni \varphi \longmapsto \Phi \circ \varphi \circ \Phi^{-1} \in \hbox{Aut}(\Omega_2)
$$
is a one-to-one, surjective group homomorphism (i.e., an {\it isomorphism}).  In particular, if two domains have
automorphism groups which are not isomorphic, then the two domains cannot be biholomorphic.
This is the strategy that we shall use to show that $B$ and $B^n$ are not biholomorphic.
For simplicity, and with no loss of generality, we restrict attention to
the case $n = 2$.  

We shall leave certain small portions of this argument as exercises for the reader, but all
the key ideas will be indicated.

Notice that, for any $\alpha, \beta \in \CC$, $|\alpha|, |\beta| < 1$, the mapping
$$
\psi_{\alpha, \beta}:  (z_1, z_2) \longmapsto \left ( \frac{z_1 - \alpha}{1 - \overline{\alpha} z_1},
                                 \frac{z_2 - \beta}{1 - \overline{\beta} z_2} \right )
$$
is a biholomorphic self-mapping (i.e., an automorphism) of the bidisc $D^2$.  We assume,
seeking a contradiction, that $\Phi: B \ra D^2$ is in fact a biholomorphic mapping
of $B$ with $D^2$.  If $(\alpha, \beta) = \Phi(0,0)$, then $\psi_{\alpha, \beta} \circ \Phi$ 
is also a biholomorphic mapping of $B$ to $D^2$ that takes the origin $(0,0)$ to the origin
$(0,0)$.  For convenience we shall denote this new mapping with the same name $\Phi$.

For any domain $\Omega$ with $0 \in \Omega$, we let $I_0^\Omega$ denote the {\it isotropy group}
of the origin:  this is the subgroup of $\hbox{Aut}(\Omega)$ consisting of those mappings
which fix the origin.  Obviously the putative biholomorphic mapping $\Phi$ of $B$ to
$D^2$ maps $I_0^B$ in $\hbox{Aut}(B)$ isomorphically onto $I_0^{D^2}$ in $\hbox{Aut}(D^2)$.  We
shall show now that this is impossible.  

In fact the connected component of the identity in each of these subgroups will also be mapped
isomorphically to each other.  So let us consider what those might be.  In the bidisc, the maps
that fix the origin are (by a simple application of the Schwarz lemma), of two kinds:
\begin{itemize}
\item Maps which are just rotations in each variable separately.
\item Maps which permute the two variables.
\end{itemize}
And of course we also must allow compositions of these two types of mappings.  But if we
are to consider the connected component of the identity in the isotropy group of the origin
then we must restrict attention to the first kind.  Call the group $G_{D^2}$.  Then we see
that
$$
G_{D^2} = \hbox{(maps which are rotations in each variable separately)} \, .
$$

Now let us look at the ball $B$.  The maps that fix the origin are of two kinds:
\begin{itemize}
\item Unitary rotations.
\item Antipodal mappings.
\end{itemize}
Again, if we wish to consider only the connected component of the identity, then
we can only consider mappings of the first kind.  Thus
$$
G_B = \hbox{(unitary rotations)} \, .
$$

Now the key fact to notice here is that $G_{D^2}$ is {\it abelian}, while
$G_B$ is {\it not abelian}.  Thus these two groups cannot be isomorphic.
We conclude then that $D^2$ and $B$ cannot be biholomorphic.  That is Poincar\'{e}'s 
theorem.

\section{Partial Differential Equations}

In the present section we use ideas from Section 2 to derive a formula
for the solution of the inhomogeneous Cauchy-Riemann equation.  This
will be a new idea for most readers.  Afterwards we shall supply an 
interesting application of the ideas.

Recall that the classical Cauchy-Riemann equations for a function
$f(z) = u(z) + iv(z)$ are
$$
\frac{\partial u}{\partial x} = \frac{\partial v}{\partial y} \qquad 
\hbox{and} \qquad \frac{\partial u}{\partial y} = - \frac{\partial v}{\partial x} \, .
$$
Expressed in terms of the complex derivative, these equations may be written neatly
as
$$
\frac{\partial}{\partial \overline{z}} f \equiv 0 \, .
$$

It turns out to be a matter of considerable interest to solve the equation
$$
\frac{\partial}{\partial \overline{z}} f \equiv \alpha 
$$
for a suitable function $\alpha$.  We shall prove the following result.

\begin{theorem}    
Let $\alpha \in C^1_c(\CC).$  The function defined
by 
$$ 
f(\zeta) = - \frac{1}{2\pi i} \int \frac{\alpha(\xi)}{\xi - \zeta} d\bar{\xi} \wedge d\xi
            = - \frac{1}{\pi} \int \frac{\alpha(\xi)}{\xi - \zeta} dA(\xi) \, ,
$$
where $dA$ is area measure, satisfies
$$  
\dbar f(\zeta) = \frac{\partial f}{\partial \bar{\zeta}}(\zeta) d\bar{\zeta} = \alpha(\zeta) d\overline{\zeta} .  
$$
\end{theorem}
\noindent {\bf Proof:}  Let $D(0,R)$ be a large disc that contains the support of
$\alpha.$  Then
\begin{eqnarray*}
  \frac{\partial f}{\partial \overline{\zeta}} (\zeta) & = & - \frac{1}{2\pi i} \frac{\partial}{\partial \overline{\zeta}}
                                             \int_\CC \frac{\alpha(\xi)}{\xi - \zeta} d\bar{\xi} \wedge d\xi  \\
                      & = & - \frac{1}{2\pi i} \frac{\partial}{\partial \overline{\zeta}}
                                             \int_\CC \frac{\alpha(\xi + \zeta)}{\xi} d\bar{\xi} \wedge d\xi  \\
                      & = & - \frac{1}{2\pi i} \int_\CC 
                               \frac{\frac{\partial\alpha}{\partial \bar{\xi}}(\xi + \zeta)}{\xi} 
                                  \d\bar{\xi} \wedge d\xi  \\
                      & = & - \frac{1}{2\pi i} \int_{D(0,R)} 
                              \frac{\frac{\partial\alpha}{\partial \bar{\xi}}(\xi)}{\xi - \zeta} 
                                  \d\bar{\xi} \wedge d\xi .
\end{eqnarray*}
By Theorem 3, this last equals
$$ 
\alpha(\zeta) - \frac{1}{2\pi i} \int_{\partial D(0,R)} \frac{\alpha(\xi)}{\xi - \zeta} d\xi = \alpha(\zeta) . 
$$
Here we have used
the support condition on $\alpha.$  This is the result that we wish to prove.  
\endpf 
\medskip \\

A simple limiting argument shows that the theorem remains valid if the
function $\alpha$ is only bounded (and not smooth with compact support).
We omit the details of the proof. This additional observation will prove
useful in the example that follows.

In general the function $f$ that we produce in the last theorem as a
solution of the $\overline{\partial}$ equation will {\it not} be compactly
supported---even thought $\alpha$ itself {\it is} compactly supported. See
[KRA3] for the details of this assertion.   It follows from elementary estimates
that the solution $f$ {\it is} bounded.  We now give a simple example of the utility
of these partial differential equations ideas.

\begin{example} \rm
Let $\Omega \subseteq \CC$ be a domain.  Let $P \in \partial \Omega$ and $U$ a small,
open neighborhood of $P$.  Suppose that $h$ is a function, holomorphic on $U \cap \Omega$,
such that $h$ blows up at $P$ at a prescribed rate.  Then there is a holomorphic function $\widehat{h}$ on
$\Omega$ that blows up at $P$ at that same prescribed rate.

To see this, let $\varphi$ be a $C^\infty_c$ function with support in $U$ that is identically equal
to 1 in a neighborhood of $P$.  Define 
$$
\alpha(z) = \left ( \frac{\partial \varphi}{\partial \overline{z}} \right ) \cdot h \, .
$$
Then certainly $\alpha$ is a bounded function---for $\partial \varphi/\partial \overline{z}$ is identically
0 in a neighborhood of $P$.  Thus the equation
$$
\frac{\partial f}{\partial \overline{z}} = \alpha
$$
has a bounded solution.  But then the
function
$$
\widehat{h} = \varphi \cdot h - f
$$
will satisfy the Cauchy-Riemann equations---since
$$
\frac{\partial}{\partial \overline{z}} \left [ \varphi \cdot h - f \right ] =
        \frac{\partial \varphi}{\partial \overline{z}} \cdot h -  \frac{\partial \varphi}{\partial \overline{z}} \cdot h \equiv 0 \, .
$$

Moreover, we see that $\widehat{h}$ will blow up at $P$ at just the same rate as $h$ itself, since the two
functions differ there by the bounded function $f$.  That completes our analysis.
\endpf
\end{example}

In this example we have been deliberately nonspecific about what ``blows up at a prescribed
rate'' means.  We are endeavoring to avoid the question of whether the singularity is
a pole or an essential singularity.  In point of fact, this is {\it not} an isolated
singularity in the usual sense of complex variable theory (see [GRK]).  Rather, this
particular singularity is at the boundary---so the familiar analyses do not apply.
										  
\section{Geometry}

Ever since Lars Ahlfors's seminal paper [AHL]
on the Schwarz lemma, geometry has played a decisive
role in complex function theory.  Today it is one of
the major tools.  In this section we show how
Poincar\'{e}'s theorem on the biholomorphic inequivalence
of the ball and the polydisc (already treated in Section 5 from
a different point of view) may be derived using
geometric ideas.  This is both a validation of the geometric
viewpoint and an introduction to some important ideas that
grow out of the Riemann mapping theorem.

As we have done in the past, we restrict attention to dimension $n=2$ for simplicity
and convenience.  The reader may check that the arguments go through unchanged in
all dimensions.

We begin now with some terminology.

\begin{definition} \rm  If  $U_{1},\ U_{2} \subseteq \CC^{2}$ 
and
$$ 
f = (f_1, f_2): U_{1} \longrightarrow U_{2} 
$$ 

\noindent is holomorphic then define  $\hbox{\rm Jac}_\CC f(z)
\ , z  \in   U_{1},$ (the {\it Jacobian matrix} of $f$ at $P$) to
be the matrix
$$ 
\left (  \matrix{ {\partial f_{1} \over {\partial z_{1}}} (z) 
& {\partial f_{1} \over {\partial z_{2}}} (z) \cr
{\partial f_{2} \over {\partial z_{1}}} (z) 
& {\partial f_{2} \over {\partial z_{2}}} (z) \cr} \right) 
$$
\end{definition}

\begin{remark} \rm
Notice that  $\hbox{\rm Jac}_ \CC f(z)$  is
distinct from the real Jacobian, $\hbox{\rm Jac}_ \RR f$, of
calculus:  the latter would be a  $4 \times 4$ matrix which arises
from treating  $f$  as a function from a domain in $\RR^{4}$ 
to a domain in  $\RR^{4}.$
\endpf
\end{remark}

We continue to use the symbol  $D$  to denote the unit disc in 
$\CC.$ Let  $\Omega  \subseteq \CC^{2}$  be a domain and  $P \in
\Omega.$ Define $(D,\Omega)_{P}$  to be the holomorphic functions  
$f: \Omega \rightarrow D$  such that  $f(P) = 0.$
Define  $(\Omega,D)_{P}$  to be the holomorphic functions  
$f: D \rightarrow \Omega$  such that  $f(0) = P.$
We are now going to define the Carath\'eodory and Kobayashi
metrics.  Since we are working in a two-dimensional space, 
we can no longer specify a metric as a scalar-valued function on
the domain.  In fact a metric will measure the length of a vector
at a point.  Here and throughout, $ |\ \ | $ denotes Euclidean
length.

\begin{definition} \rm If  $P  \in  \Omega$  and  $\xi  \in \CC^{2}$
define the {\it Carath\'{e}odory length} of $\xi$ at $P$ to be
$$
F^\Omega_C (P,\xi ) = \sup \{ |\hbox{\rm Jac}_\CC f(P) \xi|:f \in
(D,\Omega)_P \} \, . 
$$
\end{definition}

\begin{definition} \rm
If  $P  \in  \Omega$  and  $\xi  \in {\bf
C}^{2}$  define
the {\it Kobayashi length} of $\xi$ at $P$ to be

\begin{eqnarray*}
F^\Omega_K (P,\xi ) & = & \inf \{ | \xi | / | g'(0)| : g \in (\Omega,D)_{P}, \\
& &\qquad \qquad g'(0) \hbox{ is a scalar multiple of } \xi \} \, .
\end{eqnarray*}
\end{definition}

The way that we define metrics is motivated by Riemann's
philosophy.  If  $\gamma: [0,1] \rightarrow \Omega$  is a continuously
differentiable curve, we define its Kobayashi length to be
$$ 
\ell _{K} (\gamma ) = \int^1_0  F^\Omega_K (\gamma (t),\gamma' (t))
dt \, . 
$$

\noindent Notice that we are integrating the lengths, 
in the metric, of the tangent vectors to the curve.
The Carath\'eodory length of a curve is defined similarly.

One of the basic properties that we shall prove is that holomorphic
mappings decrease distance in the Carath\'eodory and Kobayashi
metrics.  In several variables we express this assertion as
follows.

\begin{proposition} \sl  Let  $U_{1}$  and  $U_{2}$  be domains and
$$ f: U_{1} \longrightarrow U_{2} $$  
be a holomorphic mapping.  If $ P  \in  U_{1}$  and  $\xi \in  \CC^{2}$ 
we define

$$ f_{*} (P) \xi  = \hbox{\rm Jac}_ \CC f(P) \xi . $$
Then

$$ F^{U_1}_C (P,\xi ) \geq   F^{U_2}_C (f(P),f_{*} (P) \xi ) $$
and

$$ F^{U_1}_K (P,\xi ) \geq  F^{U_2}_K (f(P),f_{*} (P) \xi ) .$$
\end{proposition}
\noindent {\bf Proof:}  We give the proof for the Carath\'eodory metric.
The proof for the Kobayashi metric is similar.

Choose  $\varphi \in  (D,U_{2})_{f(P)}.$
Then  $\varphi \circ  f  \in   (D,U_{1})_{P}.$
Hence

\begin{eqnarray*}
F^{U_1}_C (P,\xi ) & \geq & |  (\hbox{\rm Jac}_\CC
(\varphi \circ f) (P)) \xi | \\
& = & |  \hbox{\rm Jac}_\CC \varphi (f(P)) \circ (\hbox{\rm
Jac}_\CC f(P)) \xi | \\
& = & |  \hbox{\rm Jac}_\CC \varphi (f(P)) (f_{*} (P)) \xi | \, .
\end{eqnarray*}

\noindent Taking the supremum over all  $\varphi $  gives
$$ 
F^{U_1}_C (P,\xi ) \geq  F^{U_2}_C (f(P),f_{*} (P) \xi ) \, .  \eqno \BoxOpTwo
$$		 

\begin{corollary} \sl If  $f$  is a biholomorphic map then it 
preserves both the Carath\'eodory and the Kobayashi metrics; 
that is, the inequalities in the proposition become equalities.
\end{corollary}
\noindent {\bf Proof:}  Obvious.
\endpf 
\medskip \\

\noindent {\bf Exercise:}  Check that the proposition implies that  $f$ 
decreases the lengths of curves.  That is, if  $\gamma $  is a
continuously differentiable curve in $U_{1}$ and  $f_{*} \gamma 
\equiv f \circ \gamma$ is the corresponding curve in  $U_{2}$
then show that $\ell _{C} (f_{*} \gamma) \leq \ell _{C} (\gamma)$ 
and $\ell _{K} (f_{*} \gamma ) \leq  \ell _{K} (\gamma ).$ 
\endpf 
\medskip \\

We now define two interesting new invariants.
These invariants were trivial in one complex variable, but now they
provide crucial information.

\begin{definition} \rm
Let  $\Omega  \subseteq \CC^{2}$  be a
domain and $P  \in \Omega.$  The  {\it Carath\'eodory indicatrix}   of 
$\Omega$  at $P$  is 
$$ 
{\bf i}^C_P (\Omega) = \{ \xi \in \CC: F^\Omega_C (P,\xi ) < 1 \} \, .
$$

\noindent The  {\it Kobayashi indicatrix}   of  $\Omega$  at  $P$  is
$$ 
{\bf i}^K_P (\Omega) = \{ \xi \in \CC: F^\Omega_K (P, \xi ) < 1 \} \, .
$$

In words, the indicatrix is the ``unit ball,''\ in the indicated
metric, of vectors at  $P.$
\end{definition}

\smallskip
\begin{proposition} \sl
Let $f: \Omega_{1} \rightarrow \Omega_{2}$  
be a biholomorphic mapping of domains in  $\CC^{2}.$
Say that  $f(P) = Q.$
Then
$$ 
\hbox{\rm Jac}_ \CC f(P): {\bf i}^C_P (\Omega_{1})
\longrightarrow {\bf i}^C_Q (\Omega_{2}) 
$$  

\noindent and
$$ 
\hbox{\rm Jac}_ \CC f(P): {\bf i}^K_P (\Omega_{1})
\longrightarrow {\bf i}^K_P (\Omega_{2}) 
$$  

\noindent are linear isomorphisms.
\end{proposition}

\noindent {\bf Proof:} Since  $f$  is distance decreasing in the Kobayashi metric, 
$\hbox{\rm Jac}_ \CC f(P)$ maps ${\bf i}^K_P (\Omega_{1})$ into
${\bf i}^K_Q (\Omega_{2}).$  But the same observation applies to
$$ 
\hbox{\rm Jac}_ \CC (f^{-1}) (Q) = (\hbox{\rm Jac}_ \CC
f(P))^{-1} \, ; 
$$  

\noindent it maps ${\bf i}^K_Q (\Omega_{2})$  into  ${\bf i}^K_P
(\Omega_{1}).$  Thus  $\hbox{\rm Jac}_ \CC f(P)$  is a linear
isomorphism of ${\bf i}^K_P$ to  ${\bf i}^K_Q$ as claimed.

The proof for  ${\bf i}^C_P$  is identical.
\endpf
\medskip \\

\begin{proposition} \sl
Let  $B = B(0,1)$ be the unit ball.
Then we have
$$ 
{\bf i}^K_0 (B) = B \, .
$$
\end{proposition}

\noindent {\bf Proof:}  Let  $\varphi  \in   (B,D)_{0}.$
If  $\eta $  is any Euclidean unit vector in $\CC^{2}$  then
consider the function
$$ 
h(\zeta ) \equiv  \varphi  (\zeta ) \cdot \eta \, . 
$$

\noindent Here ``$ \cdot $'' denotes the usual inner product of
2-vectors.  We have that  $h$  maps the disc to the disc and  $h(0)
= 0.$  By the Schwarz lemma of one variable,
$$ 
|  h'(0) | \leq 1 \, . 
$$

\noindent Since this inequality holds for any choice of  $\eta,$
we conclude that
$$ 
|  \varphi' (0) |   \leq   1 \, . 
$$

Now if  $\xi $  is any vector in  $\CC^{2}$  then it follows
from the
preceding calculation that

\begin{eqnarray*} 
F^B_K (0,\xi ) 
& = & \inf \{ |  \xi |  / |  \varphi' (0) | : \varphi \in  (B,D)_{0}
\} \\  
& \geq &  |  \xi | \, . 
\end{eqnarray*}

\noindent On the other hand, the map
$$ 
\varphi_{0} (\zeta ) \equiv {\zeta \over{|\xi|}}  \xi  
$$  

\noindent satisfies  $\varphi_{0} \in  (B,D)_{0}$  and  
$\varphi_{0}' (0)$  is a positive multiple of  $\xi .$
Therefore
$$ 
F^B_K (0,\xi ) \leq   |  \xi |  / |  \varphi_{0}' (0) |  = | 
\xi | \, . 
$$

\noindent We conclude that
$$ 
F^\Omega_K (0,\xi ) = |  \xi | \, , 
$$  

\noindent hence that
$$ 
{\bf i}^K_0 (B) = B \, .  \eqno \BoxOpTwo
$$

\begin{proposition}  \sl
Let  $D^2 = D^{2} (0,1)$ be the unit bidisc.
Then
$$ 
{\bf i}^K_0 (D^2) = D^2 \, . 
$$
\end{proposition}

\noindent {\bf Proof:}  Define the projections
$$ 
\pi_1 (z_1, z_2) = z_1 \qquad \hbox{and}\qquad \pi_2 (z_1 , z_2)
= z_2 \, . 
$$

\noindent Let  $\eta  = (\eta_{1},\eta_{2}) \in  \CC^{2}$  be
any vector.  By Proposition 3, we have that
$$ 
F^{D^2}_K (0,\eta ) \geq 
F^{\pi_{1}(D^2)}_{K} (\pi_{1}(0),(\pi_{1})_{*} \eta ) = 
F^D_K (0,\eta_{1}) \, . 
$$

\noindent But the Schwarz lemma of one variable tells us easily
that the last quantity is just  $ | \eta_{1} |.$ 
A similar argument shows that
$$ 
F^{D^2}_K (0,\eta )  \geq   | \eta_{2} | \, . 
$$

\noindent We conclude from these two inequalities that
$$ 
F^{D^2}_K (0,\eta )  \geq   \max \{ | \eta_{1} | , | \eta_{2}
|  \} \, .
$$

\noindent Therefore
$$ 
{\bf i}^K_0 (D^{2}) \subseteq   D^{2} \, . 
$$

For the reverse inclusion, fix  $\eta $  as above and consider the
function
$$ 
\varphi (\zeta ) = 
\Biggl({\zeta \eta_{1} \over {\max \{ | \eta_{1} |, | \eta_{2} |
\}}},\  
{\zeta \eta_{2} \over {\max \{ | \eta_{1} |, | \eta_{2} | \}} }
\Biggr) \, . 
$$

\noindent Then it is obvious that  $\varphi \in (\Omega,D)_{0}$  
and that  $\varphi' (0)$  is a positive multiple of  $\eta.$
Therefore
$$ 
F^{\Omega}_K (0,\eta ) \leq { | \eta | \over {|\varphi'(0)|}} 
= \max \{ | \eta_{1} | , |  \eta_{2} | \} \, . 
$$

\noindent The opposite inclusion now follows.
\endpf 
\medskip \\

\noindent {\bf Exercise:} 
Verify that ${\bf i}^C_0 (B) = B$  and 
${\bf i}^C_0 (D^{2}) = D^{2}.$ 
\endpf
\medskip \\

\begin{theorem}[Poincar\'{e}]  There is no biholomorphic mapping
of the bidisc  $D^{2}$  to the ball $B.$
\end{theorem}

\noindent {\bf Proof:}  Suppose that
$$ 
\Phi : D^{2} \longrightarrow B 
$$  

\noindent is biholomorphic.  Let  $\Phi^{-1} (0) = \alpha  \in  
D^{2}.$ There is an element  $\varphi \in  \hbox{\rm Aut}( D^{2}
(0,1))$  such that  $\varphi (0) = \alpha.$ Consider  $g \equiv \Phi
\circ \varphi.$  Then
$$ g: D^{2} \longrightarrow  B 
$$  

\noindent is biholomorphic and  $g(0) = 0.$  We will show that  $g$
cannot exist.

By Proposition 11, $\hbox{\rm Jac}_ \CC g(0)$  is a linear
isomorphism of ${\bf i}^K_0 (D^{2})$  to   ${\bf i}^K_0 (B).$
But Propositions 12 and 13 identify these as  $D^{2}$  and  $B$ 
respectively.  So we have that
$$ 
\hbox{\rm Jac}_ \CC g(0): D^{2} \longrightarrow  B 
$$  

\noindent is a linear isomorphism.  However this is impossible.
For the segment  ${\bf i} = \{ (t + i0,1): 0  \leq   t \leq   1 \}$
lies in  $\partial D^{2}.$ The linear isomorphism would map ${\bf
i}$   to a nontrivial segment in  $\partial B.$  But  $B$  is
strictly convex (all boundary points are extreme) so its boundary
contains no segments.  This is the required contradiction.
\endpf
\medskip \\

\begin{remark}\rm  It is important to appreciate the logic in this
proof.  The hypothesized map  $\Phi $ (and therefore  $g$  as well)
is {\it not }   assumed to extend in any way to  $\partial D^{2}.$
Indeed, given the very different natures of  $\partial D^{2}$  and 
$\partial B,$ we would expect  $\Phi $  to be highly pathological
at the boundary.  Our geometric machinery allows us to pass to the
linear map  $\hbox{\rm Jac}_ \CC g(0),$ which is defined on all
of space.  Thus we are able to analyze the boundaries of the
domains and arrive at a contradiction.
\endpf
\end{remark}

\section{Concluding Remarks}

We have given several examples to show how a problem that began as a complex
analysis problem could ultimately be resolved by techniques from another part
of mathematics.  This is not to say that one could not have brought
function-theoretic techniques to bear on the original problems.  But it is
always fruitful to bring in new ideas and new perspectives.

In fact it is safe to say that the application of techniques from other
disciplines has helped to spawn entirely new directions in mathematics.
Today the subjects of K\"{a}hler geometry, the $\overline{\partial}$-Neumann
problem, and automorphism groups---to give just a few examples---are 
branches of complex analysis that have a life of their own.

Moreover, those new ideas will often bear additional benefits.  They will
give new insights, offer different points of view, or (in the best
of all possible circumstances) suggest new problems.  Certainly
that is the way that matters have played out in practice in
modern complex function theory.  

Other parts of mathematics have seen similar benefits from bringing in ideas
from the outside.  In recent years we have seen applications of techniques
from partial differential equations in algebraic geometry, applications
of Lie group theory in several complex variables, and applications
of abstract logic in real analysis.  The number of examples keeps growing,
and mathematics is correspondingly enriched.  We hope that this
article has served as a modest introduction to this exciting
set of ideas.

\newpage

\noindent {\Large \sc References}
\smallskip \\

\begin{enumerate}

\item[{\bf [AHL]}]  L. Ahlfors, An extension of Schwarz's lemma
{\it Trans. Amer. Math. Soc.} 43 (1938), 359--364.

\item[{\bf [BLK]}]  B. Blank and S. G. Krantz, {\it Calculus, Multivariable}, Key College Press,
Emeryville, CA, 2006.

\item[{\bf [GIT]}]  Gilbarg and N. Trudinger, {\it Elliptic Partial Differential Equations of Second Order},
$2^{\rm nd}$ ed., Springer-Verlag, Berlin, 2001.

\item[{\bf [GRK]}]  R. E. Greene and S. G. Krantz, {\it Function Theory of One
Complex Variable}, $3^{\rm rd}$ ed., American Mathematical Society,
Providence, RI, 2006.

\item[{\bf [HOR]}]  L. H\"{o}rmander, {\it Linear Partial Differential Equations}, Springer-Verlag,
Berlin, 1969.  

\item[{\bf [KOB]}]  S. Kobayashi, {\it Hyperbolic Manifolds and Holomorphic Mappings},
Marcel Dekker, New York, 1970.

\item[{\bf [KRA1]}]  S. G. Krantz, The Carath\'{e}odory and Kobayashi metrics 
and applications in complex analysis, {\it American Mathematical Monthly}, to appear.
							
\item[{\bf [KRA2]}]  S. G. Krantz, {\it Partial Differential Equations and Complex Analysis},
CRC Press, Boca Raton, FL, 1992.

\item[{\bf [KRA3]}]  S. G. Krantz, {\it Function Theory of Several Complex Variables},
$2^{\rm nd}$ ed., American Mathematical Society, Providence, RI, 2001.

\item[{\bf [KRA4]}]  S. G. Krantz, {\it Geometric Function Theory:  Exploration in Complex Analysis},
Birkh\"{a}user, Boston, 2006.

\item[{\bf [KRA5]}]  S. G. Krantz, {\it A Panorama of Harmonic Analysis}, Mathematical Association of
America, Washington, D.C., 1999.

\item[{\bf [SIC]}] J. Siciak, On series of homogeneous polynomials and
their partial sums, {\it Annales Polonici Math.} 51(1990), 289--302.

\item[{\bf [STE]}]  E. M. Stein, {\it Singular Integrals and Differentiability Properties
of Functions}, Princeton University Press, Princeton, NJ, 1970.

\end{enumerate}
\vspace*{.25in}

\noindent {\bf STEVEN G. KRANTZ} received his B.A. degree from
the University of California at Santa Cruz in 1971. He earned
the Ph.D. from Princeton University in 1974. He has taught at
UCLA, Princeton University, Penn State, and Washington
University in St.\ Louis. Krantz is the holder of the UCLA
Alumni Foundation Distinguished Teaching Award, the Chauvenet
Prize, and the Beckenbach Book Prize. He is the author of 150
papers and 50 books. His research interests include complex
analysis, real analysis, harmonic analysis, and partial
differential equations. Krantz is currently the Deputy
Director of the American Institute of Mathematics.
\medskip \\
{\it American Institute of Mathematics, 360 Portage Avenue, \\
Palo Alto, CA 94306} 
\medskip \\
{\it sk@math.wustl.edu}

\end{document}